\documentclass[12pt, a4paper]{article}
\usepackage{latexsym}
\usepackage{amsmath}
\usepackage{amssymb}
\usepackage{color}
\usepackage{amsfonts}
\usepackage{verbatim}
\usepackage{amsthm}

\usepackage{anysize,hyperref}
\input xypic
\xyoption{all}

\usepackage[perpage,symbol]{footmisc}
\usepackage{setspace}
\renewcommand{\cal}{\mathcal}

\newcommand{\rtri}[6]{\xymatrix@C=1.5em{#1\ar[r]^{#4}&#2\ar[r]^{#5}&#3\ar[r]^{#6}&\Sigma #1}}
\newcommand{\tri}[7]{\xymatrix@C=1.5em{#1\ar[r]^{#5}&#2\ar[r]^{#6}&#3\ar[r]^{#7}&#4}}

\begin{document}
\baselineskip=15pt
\title{\large{\bf Subfactor categories of triangulated categories}}
\medskip
\author{Jinde Xu \qquad Panyue Zhou\qquad Baiyu Ouyang \footnotemark[1]{}\\
{\small College of Mathematics and Computer Science},\\
{\small Key Laboratory of High Performance Computing and}\\
{\small Stochastic Information Processing (Ministry of Education of China)},\\
{\small Hunan Normal University, Changsha, Hunan 410081, P. R.
China}\\}

\date{}

\maketitle
\def\blue{\color{blue}}
\def\red{\color{red}}

\newtheorem{theorem}{Theorem}[section]
\newtheorem{lemma}[theorem]{Lemma}
\newtheorem{corollary}[theorem]{Corollary}

\newtheorem{proposition}[theorem]{Proposition}

\newtheorem{conjecture}{Conjecture}

\theoremstyle{definition}
\newtheorem{definition}[theorem]{Definition}
\newtheorem{question}[theorem]{Question}
\newtheorem{remark}[theorem]{Remark}
\newtheorem{remark*}[]{Remark}
\newtheorem{example}[theorem]{Example}
\newtheorem{example*}[]{Example}

\baselineskip=15pt
\parindent=0.5cm

\footnotetext[1]{\noindent Corresponding author.\\
\hspace*{2em}Email Address: hnu\_xujinde@126.com(J.D.Xu),
 panyuezhou@163.com(P.Y.Zhou),\\oy@hunnu.edu.cn(B.Y.Ouyang).\\
\hspace*{2em}Supported by the NSF of China (No. 11371131) and
Construct Program of the Key Discipline in Hunan Province.}

\begin{abstract}
\baselineskip=16pt Let $\cal T$ be a triangulated category, $\cal A$
 a full subcategory of $\cal T$ and $\cal X$ a
functorially finite subcategory of $\cal A$. If $\cal A$ has the
properties that any $\cal X$-monomorphism of $\cal A$ has a cone and
any $\cal X$-epimorphism has a cocone. Then the subfactor category
$\cal A/[\cal X]$ admits a pretriangulated structure in the sense of
\cite{br}. Moreover the above pretriangulated category $\cal A/[\cal
X]$ with $(\cal X, \cal X[1])=0$ becomes a triangulated category if
and only if $(\cal A,\cal A)$ forms an $\cal X$-mutation pair and
$\cal A$ is closed under extensions.
\bigskip

\medskip
\end{abstract}

\section{Introduction}
\indent Over the past decades, triangulated categories have made
their way into many different parts of mathematics and have become
indispensable in many different areas of mathematics. Nowadays there
are important applications of triangulated categories in areas like
algebraic geometry, algebraic topology, commutative algebra,
differential geometry, microlocal analysis or representation theory.

The most influenced work of triangulated categories in
representation theory was created by Happel \cite{h}, he proved that
the stable category of a Frobenius category is a triangulated
category in last two decades. Later Beligiannis, Marmaridis and Reiten et al.
studied the one-side triangulated categories in a series of works
\cite{abm,bm,br,b}, one of the importance results of them is that
any contravariantly (resp. covariantly) finite subcategory $\cal X$
of $mod\Lambda$ induces on the stable category
$\underline{mod}_{\cal X}\Lambda$ (resp. $\overline{mod}_{\cal
X}\Lambda$) of $mod\Lambda$ a left (resp. right) triangulated
category, where $mod\Lambda$ is the category of finitely generated
$\Lambda$-modules over an artin algebra $\Lambda$. In the recent
paper \cite{iy}, Iyama and Yoshino proved that if $\cal D\subset
\cal Z$ are extension closed subcategories of a triangulated
category $\cal T$ with $\cal D$ satisfying $(\cal D,\cal D[1])$=0,
and if $(\cal Z,\cal Z)$ is a $\cal D$-mutation pair. Then \emph{the
subfactor category} $\cal Z/[\cal D]$ called by Iyama and Yoshino is
a triangulated category. Soon later, J{\o}rgensen \cite{j} gave a
similar construction of triangulated category by quotient category
in another manner, and Liu and Zhu\cite{lz} applied these
constructions to the one-side triangulated categories.

In this paper, we construct a one-side triangulated structure for the subfactor category
$\cal T/[\cal X]$ where $\cal X$ is a contravariantly or covariantly
finite subcategory of a triangulated category $\cal T$.

\noindent\textbf{Main Theorem 2.9.}

Let $\cal A$ be a full subcategory of a triangulated category $\cal
T$, $\cal X$ a covariantly finite subcategory and $\cal Y$ a
contravariantly finite subcategory of $\cal A$
\begin{itemize}
\item[1.] If the subcategory $\cal A$ has the property that any $\cal X$-monomorphism of $\cal A$
has a cone, then the subfactor category $\cal A/[\cal X]$ forms a right
triangulated category.
\item[2.] If the subcategory $\cal A$ has the property that any $\cal Y$-epimorphism of $\cal A$
has a cocone, then the subfactor category $\cal A/[\cal Y]$ forms a left
triangulated category.
\end{itemize}

Moreover, If $\cal X$ is a functorially finite subcategory of $\cal
A$ and $\cal A$ has the properties that any $\cal X$-monomorphism of
$\cal A$ has a cone and any $\cal X$-epimorphism has a cocone, then
the subfactor category $\cal A/[\cal X]$ forms a pretriangulated
category. Noting that when the subcategory $\cal A$ is $\cal T$, the
theorem is just the Theorem 1.2 of \cite{j}. Later, we apply our
results to the $\cal D$-mutation setting of \cite{iy}, and get the
following result: The above pretriangulated category $\cal A/[\cal
X]$ with $(\cal X, \cal X[1])=0$ becomes a triangulated category if
and only if $(\cal A,\cal A)$ forms an $\cal X$-mutation pair and
$\cal A$ is closed under extensions.

\bigskip

\section{Subfactor categories} Throughout this paper we assume,
unless other stated, that all consider categories are $k$-linear
Hom-finite, skeletally small, and Krull-Schmidt, where $k$ is a
field. We denote by $Hom_{\cal C}(X,Y)$ or $\cal C(X,Y)$ the set of
morphisms from $X\to Y$ in a category $\cal C$. When we say that
$\cal D$ is a subcategory of $\cal C$, we always mean that $\cal D$
is a full subcategory which is closed under isomorphisms, direct
sums and direct summands.

We begin by recall some definitions and notations of approximations
and homologically finite subcategory of an arbitrary
category\cite{ar}. More information please refer to \cite{ar}.

Let $\cal A$ be a category and $\cal X$ a subcategory of $\cal A$. A
morphism $f:X_B\to B$ of $\cal A$ with $X_B$ an object in $\cal X$,
is said to be a \emph{right $\cal X$-approximation of $B$}, if the
morphism $\cal A(X,f_B):\cal A(X,X_B)\to\cal A(X,B)$ is surjective
for all objects $X$ in $\cal X$. The subcategory $\cal X$ is said to
be a \emph{contravariantly finite subcategory of $\cal A$} if any
object $B$ of $\cal A$ has a right $\cal X$-approximation. Dually a
\emph{left $\cal X$-approximation} and a \emph{covariantly finite
subcategory} of $\cal A$ are defined. A contravariantly and
covariantly finite subcategory is called \emph{functorially finite}.

Given two objects $A$ and $B$ of $\cal A$, we denote by $[\cal
X](A,B)$ the set of morphisms from $A$ to $B$ of $\cal A$ which
factor through some object of $\cal X$. It is well known that $[\cal
X](A,B)$ is a subgroup of $\cal A(A,B)$, and that the family of
these subgroups $[\cal X](A,B)$ forms an ideal of $\cal A$. Thus we
have the category $\cal A/[\cal X]$ whose objects are objects of
$\cal A$ and whose morphisms are elements of $\cal A(A,B)/[\cal
X](A,B)$. The composition of $\cal A/[\cal X]$ is induced
canonically by the composition of $\cal A$. We denote by $\bar{A}$
the image of an object A of $\cal A$ in $\cal A/[\cal X]$ and
$\bar{f}$ the image of $f:A\to B$ of $\cal A$ in $\cal A/[\cal
X]$.\\

For basic references on representation theory of triangulated
categories, we refer to \cite{h}.

We recall some basic notions on one-sided triangulated categories from
\cite{abm,bm,br}. Let $\cal C$ be an additive category and
$\Sigma:\cal C\to \cal C$ an additive endofunctor called
\emph{suspension functor}. A \emph{sextuple} $(A,B,C,f,g,h)$ in
$\cal C$ is given by the form of $\rtri{A}{B}{C}{f}{g}{h}$ with
$A,B,C\in \cal C$. \emph{A morphism from sextuples} $(A,B,C, f, g,
h)$ to $(A',B',C', f', g', h')$ is a triple $(\alpha,\beta,\gamma)$
of morphisms of $\cal C$, which makes the next diagram commutative:
$$\xymatrix{A\ar[r]^f\ar[d]^{\alpha}&B\ar[r]^g\ar[d]^{\beta}&C\ar[r]^h\ar[d]^{\gamma}&\Sigma A\ar[d]^{\Sigma \alpha}\\
A'\ar[r]^{f'}&B'\ar[r]^{g'}&C'\ar[r]^{h'}&\Sigma A'}$$ If in
addition $\alpha,\beta,$ and $\gamma$ are isomorphisms in $\cal C$,
the morphism $(\alpha,\beta,\gamma)$ is then called an
\emph{isomorphism of sextuples.}

The composition of the morphisms of sextuples is induced in the
canonical way by the corresponding composition of the morphisms of
$\cal C$.

\begin{definition}
A set $\nabla$ of sextuples in $\cal C$ is called \emph{a right
triangulation of $\cal C$} if it is closed under isomorphisms and
satisfies the following axioms. The elements of $\nabla$ are then
called \emph{right triangles}.

[rTR0] For any object $A$ of $\cal C$, the sextuple
$\xymatrix{0\ar[r]^0&A\ar[r]^{1_A}&A\ar[r]^0&0}$ belongs to $\nabla$

[rTR1] Every morphism $f:A\to B$ in $\cal C$ can be embedded into a
right triangle $\rtri{A}{B}{C}{f}{g}{h}$.

[rTR2] If $\rtri{A}{B}{C}{f}{g}{h}$ is a right triangle, then
$\rtri{B}{C}{\Sigma A}{g}{h}{-\Sigma f}$ is a right triangle.

[rTR3] Given two right triangles $\rtri{A}{B}{C}{f}{g}{h}$ and
$\rtri{A'}{B'}{C'}{f'}{g'}{h'}$ and morphisms $\alpha:A\to A',\
\beta:B\to B'$ such that $\beta f=f'\alpha$, there exists a morphism
$(\alpha,\beta,\gamma)$ from the first triangle to the second.

[rTR4] Given right triangles $\rtri{A}{B}{C}{f}{g}{h}$,
$\rtri{B}{X}{Y}{a}{b}{c}$ and $\rtri{A}{X}{Z}{af}{d}{e}$. Then there
exist morphism $s:C\to Z$ and $t:Z\to Y$ such that the following
diagrams commute and the third column in first diagram is a right
triangle.

$$\xymatrix{A\ar[r]^f\ar@{=}[d]&B\ar[r]^g\ar[d]^a&C\ar[r]^h\ar[d]^s&\Sigma A\ar@{=}[d]\\
A\ar[r]^{af}&X\ar[r]^d\ar[d]^b&Z\ar[r]^e\ar[d]^t&\Sigma A\\
&Y\ar@{=}[r]\ar[d]^c&Y\ar[d]\\
&\Sigma B\ar[r]&\Sigma C}$$
$$\xymatrix{A\ar[r]^{af}\ar[d]^f&X\ar[r]^d\ar[d]^1&Z\ar[r]^e\ar[d]^t&\Sigma A\ar[d]^{\Sigma f}\\
B\ar[r]^a&X\ar[r]^b&Y\ar[r]^c&\Sigma B}$$

\end{definition}
The additive category $\cal C$ together with the suspension functor
$\Sigma$ and the right triangulation $\nabla$ is called a
\emph{right triangulated category}, denoted by the triple $(\cal
C,\Sigma, \nabla)$.

The left triangulated category $(\cal C,\Omega, \triangle)$ is
defined dually, where the endofunctor $\Omega:\cal C\to \cal C$ is
called \emph{loop functor} and $\triangle$ is the set of left
triangles satisfying axioms analogous to rTR0-rTR4.\\

Let $\cal T$ be a triangulated category and $\cal X$ a subcategory
of $\cal T$. We recall that a morphism $f:A\to B$ is called
\emph{$\cal X$-monic}, if the induced morphism $\cal T(f,X):\cal
T(B,X)\to\cal T(A,X)$ is surjective for any object $X$ of $\cal X$.
Dually a \emph{$\cal X$-epimorphism} morphism is defined. Obviously,
any left $\cal X$-approximation is $\cal X$-monic and any right
$\cal X$-approximation is $\cal X$-epic.

\begin{lemma}\label{l1}
Let $l:A\to X_A$ be a $\cal X$-monomorphism and $f:A\to B$ any
morphism in $\cal T$. Then

(1)The morphism $\binom{f}{l}:A\to B\oplus X_A$ is also $\cal
X$-monic.

(2)The morphism $g$ in the following commutative diagram with
triangles for rows is also $\cal X$-monic.
$$\xymatrix@C=3em{A'[-1]\ar[r]\ar@{=}[d]&A\ar[r]^l\ar[d]^f&X_A\ar[r]&A'\ar@{=}[d]\\
A'[-1]\ar[r]&B\ar[r]^{g}&C\ar[r]&A'}$$
\end{lemma}

\proof (1) It is trivial.

(2) By Lemma 1.4.3 in \cite{n}, the diagram may be completed to a
morphism of triangles
$$\xymatrix{A'[-1]\ar[r]\ar@{=}[d]&A\ar[r]^l\ar[d]^f&X_A\ar[r]\ar@{-->}[d]^k&A'\ar@{=}[d]\\
A'[-1]\ar[r]&B\ar[r]^{g}&C\ar[r]&A'}$$ such that
$$\xymatrix{A\ar[r]^l\ar[d]^f&X_A\ar@{-->}[d]^k\\
B\ar[r]^{g}&C}$$ is homotopy cartesian. Since $l$ is $\cal X$-monic,
for any $\varphi:B\to X,\ \forall X\in \cal X$ there exists a
morphism $\psi:X_A\to X$ such that $\psi l=\varphi f$. Hence there
exists a morphism $\sigma:C\to X$ such that $\varphi=\sigma g$ by
the property of homotopy cartesian, i.e.
$g$ is $\cal X$-monic.  \qed \\

For the rest of the paper, we shall deal only with the right case,
leaving for the reader to state and prove dual results for the left
case. Let $\cal T$ be a triangulated category and $\cal X\subset
\cal A$ full subcategories of $\cal T$. From now on, we assume,
unless otherwise stated, that $\cal A$ and $\cal X$ satisfy the
following two conditions:
\begin{itemize}
\item[(A1)] $\cal X$ is a covariantly finite subcategory of $A$.
\item[(A2)] Any $\cal
X$-monomorphism of $\cal A$ has a cone. i.e. for any $\cal
X$-monomorphism $f:A\to B$ of $\cal A$, the third term $C_f$ in the
triangle $\tri{A}{B}{C_f}{A[1]}{f}{}{}$ of $\cal T$ is also in $\cal
A$.
\end{itemize}

Under such a setting, we will enrich the \emph{subfactor category}
$\cal A/[\cal X]$ of $\cal A$ a right triangulated structure.

In order to construct the right triangulation $\nabla$ of $\cal
A/[\cal X]$, First of all, we construct the \emph{suspension
functor} $\Sigma: \cal A/[\cal X]\to \cal A/[\cal X]$ as follows:
For any object $A\in \cal A$, consider the triangle in $\cal T$
$$\xymatrix{A\ar[r]^{\alpha_A}&X_A\ar[r]&\Sigma A\ar[r]&A[1]}$$
where $\alpha_A$ is a left $\cal X$-approximation of $A$ and define
$\Sigma \bar{A}$ like this. For any morphism $f\in (A,A')$, there
exist $g$ and $h$ which make the following diagram commutative.
$$\xymatrix{A\ar[r]^{\alpha_A}\ar[d]^f&X_A\ar[r]\ar[d]^g&\Sigma A\ar[r]\ar[d]^h&A[1]\ar[d]^{f[1]}\\
A'\ar[r]^{\alpha_{A'}}&X_A'\ar[r]&\Sigma A'\ar[r]&A'[1]}$$ Now put
$\Sigma \bar{f}:=\bar{h}$ and it is easy to see that $\bar{h}$ is
uniquely determined by $\bar{f}$, i.e. the endofunctor $\Sigma:\cal
A/[\cal X]\to \cal A/[\cal X]$ is well defined.

Next we construct two kinds of right triangles in $\cal A/[\cal X]$,
the distinguished and the induced ones.

The distinguished right triangles are obtained as follows:

Given a morphism $f:A\to B$ in $\cal A$, consider the following
diagram with triangles for rows in $\cal T$
$$\xymatrix{\Sigma A[-1]\ar[r]\ar@{=}[d]&A\ar[r]^{\alpha_A}\ar[d]^f&X_A\ar[r]&\Sigma A\ar@{=}[d]\\
\Sigma A[-1]\ar[r]&B\ar[r]^{g}&C_f\ar[r]^h&\Sigma A}$$ where
$\alpha_A$ is a left $\cal X$-approximation of $A$. Noting that by
Lemma \ref{l1} the morphism $\binom{f}{\alpha_A}$ is $\cal X$-monic
and
there exists a morphism $k$ such that $$\xymatrix{A\ar[r]^{\alpha_A}\ar[d]^f&X_A\ar@{-->}[d]^k\\
B\ar[r]^{g}&C_f}$$ is homotopy cartesian, i.e.
$\xymatrix{A\ar[r]^{\binom{f}{\alpha_A}\quad}&B\oplus
X_A\ar[r]^{\quad(-g,k)}&C_f\ar[r]&A[1]}$ is a triangle in $\cal T$,
hence $C_f$ is in $\cal A$.

\begin{definition}
A sextuple $\tri{\overline{X}}{\overline{Y}}{\overline{Z}}{\Sigma
\overline{X}}{\overline{u}}{\overline{v}}{\overline{w}}$ in $\cal
A/[\cal X]$ is said to be an \emph{$\cal X$-distinguished} right
triangle, if it is isomorphic to a sextuple
$\tri{\overline{A}}{\overline{B}}{\overline{C_f}}{\Sigma
\overline{A}}{\overline{f}}{\overline{g}}{\overline{h}}$ given by
some $f:A\to B$ in $\cal A$.
\end{definition}

We shall call the $\cal X$-distinguished right triangle just
distinguished right triangle whenever it is clear from the context.

\begin{remark}\label{r1}
Noting that the sextuple
$\tri{\overline{A}}{\overline{B}}{\overline{C_f}}{\Sigma
\overline{A}}{\overline{f}}{\overline{g}}{\overline{h}}$ is
independent of choice of $k$ and is uniquely determined (up to
isomorphism) by $f$. We can assume, without loss of generality, that
the sextuple is given by the following commutative diagram
$$\xymatrix{A\ar[r]^{\alpha_A}\ar[d]^f&X_A\ar[r]\ar[d]^{k}&\Sigma A\ar@{=}[d]\\
B\ar[r]^{g}&C_f\ar[r]^h&\Sigma A}$$ with triangles for rows and the
left square being homotopy cartesian.
\end{remark}

The induced right triangles are obtained as follows:

Consider the triangle $\tri{A}{B}{C}{A[1]}{f}{g}{h'}$ in $\cal T$
with $f\in \cal A$ being $\cal X$-monic, there is a commutative
diagram of triangles:
$$\xymatrix{A\ar[r]^f\ar@{=}[d]&B\ar[r]^g\ar[d]&C\ar[r]^{h'}\ar[d]^h&A[1]\ar@{=}[d]\\
A\ar[r]^{\alpha_A}&X_A\ar[r]^{\beta_A}&\Sigma
A\ar[r]^{\gamma_A}&A[1]}$$ where $\alpha_A$ is a left $\cal
X$-approximation of $A$.

\begin{remark}\label{r2}
Noting that the morphism $\overline{h}$ in $\cal A/[\cal X]$ is
unique determined by the triangle $\tri{A}{B}{C}{A[1]}{f}{g}{h'}$.
Indeed, if $a:C\to \Sigma A$ satisfies $h'=\gamma l$, then $\gamma
h=\gamma l$, i.e. $\gamma(h-l)=0$. Hence $h-l$ factor through $X_A$,
i.e. $\bar{h}=\bar{l}$.
\end{remark}

\begin{definition}
A sextuple $\tri{\overline{X}}{\overline{Y}}{\overline{Z}}{\Sigma
\overline{X}}{\overline{u}}{\overline{v}}{\overline{w}}$ in $\cal
A/[\cal X]$ is said to be an \emph{$\cal X$-induced} right triangle,
if it is isomorphic to a sextuple
$\tri{\overline{A}}{\overline{B}}{\overline{C}}{\Sigma
\overline{A}}{\overline{f}}{\overline{g}}{-\overline{h}}$ for some
$\cal X$-monomorphism $f:A\to B$ of $\cal A$.
\end{definition}

We shall call the $\cal X$-induced right triangle just induced right
triangle whenever it is clear from the context.

In fact, the set $\nabla_d$ of right distinguished triangles and the
set $\nabla_i$ of right induced triangles are equal.

\begin{proposition}\label{p1}
Any distinguished right triangle is isomorphic to an induced one and
any induced right triangle is isomorphic to a distinguished one in
$\cal A/[\cal X]$.
\end{proposition}

\proof Given a morphism $f:A\to B$ in $\cal A$, by Remark \ref{r1},
the right distinguished triangle
$\tri{\overline{A}}{\overline{B}}{\overline{C_f}}{\Sigma
\overline{A}}{\overline{f}}{\overline{g}}{\overline{h}}$ is given by
the following commutative diagram of triangles
$$\xymatrix{A\ar[r]^{\alpha_A}\ar[d]^f&X_A\ar[r]^{\beta_A}\ar[d]^{k}&\Sigma A\ar@{=}[d]\ar[r]^{\gamma_A}&A[1]\ar[d]^{f[1]}\\
B\ar[r]^{g}&C_f\ar[r]^h&\Sigma A\ar[r]&A[1]}$$ with the left square
being homotopy cartesian, that is to say there is a triangle of
$\cal T$
$$\xymatrix@C=3em{A\ar[r]^{\binom{f}{\alpha_A}\quad}&B\oplus
X_A\ar[r]^{\quad(-g,k)}&C_f\ar[r]^{\gamma_A h}&A[1]}$$ with
$A,B\oplus X_A,C_f\in \cal A$. Since $\binom{f}{\alpha_A}$ is $\cal
X$-monic, the triangle induces a right induced triangle in $\cal
A/[\cal X]$ as follows:
$$\xymatrix{A\ar@{=}[d]\ar[r]^{\binom{f}{\alpha_A}\quad}&B\oplus
X_A\ar[d]^{(0,1)}\ar[r]^{\quad(-g,k)}&C_f\ar[d]^h\ar[r]^{\gamma_A h}&A[1]\ar@{=}[d]\\
A\ar[r]^{\alpha_A}&X_A\ar[r]^{\beta_A}&\Sigma
A\ar[r]^{\gamma_A}&A[1]}.$$ Hence the induced triangle is
$\xymatrix{\bar{A}\ar[r]^{\bar{f}}&\bar{B}\ar[r]^{-\bar{g}}&\bar{C}_f\ar[r]^{-\bar{h}}&\Sigma
\bar{A}}$, obviously isomorphic to the distinguished one
$\xymatrix{\bar{A}\ar[r]^{\bar{f}}&\bar{B}\ar[r]^{\bar{g}}&\bar{C}_f\ar[r]^{\bar{h}}&\Sigma
\bar{A}}$.

Given a right induced triangle
$\xymatrix@C=1.5em{\bar{A}\ar[r]^{\bar{f}}&\bar{B}\ar[r]^{\bar{g}}&\bar{C}\ar[r]^{-\bar{\theta}}&\Sigma
\bar{A}}$, i.e. we have the following commutative diagram of
triangles
$$\xymatrix{A\ar[r]^f\ar@{=}[d]&B\ar[r]^g\ar[d]^\delta&C\ar[r]^{h}\ar[d]^{\theta}&A[1]\ar@{=}[d]\\
A\ar[r]^{\alpha_A}&X_A\ar[r]^{\beta_A}&\Sigma
A\ar[r]^{\gamma_A}&A[1]}$$ with $f, \alpha_A\in \cal A$ being $\cal
X$-monic. By Remark \ref{r2}, we can assume without loss of
generality that $\theta:C\to \Sigma$ is the morphism such that the
middle square is homotopy cartesian, then we have a triangle
$$\xymatrix@C=3em{B\ar[r]^{\binom{g}{\delta}\quad}&C\oplus
X_A\ar[r]^{\quad(\theta,-\beta_A)}&\Sigma A\ar[r]^{-f[1]\gamma_A
}&A[1]}$$ An easy computation allow us to have the following
commutative diagram of triangles.
$$\xymatrix@C=3em{A\ar[d]^{-f}\ar[r]^{\alpha_A}&X_A\ar[d]^{\binom{0}{-1}}\ar[r]^{\beta_A}&\Sigma A\ar@{=}[d]\ar[r]^{\gamma_A}&A[1]\ar[d]^{-f[1]}\\
B\ar[r]^{\binom{g}{\delta}}&C\oplus
X_A\ar[r]^{(\theta,-\beta)}&\Sigma A\ar[r]^{-f[1]\gamma_A}&B[1]}$$
Then by the definition of right distinguished triangles, the
sextuple
$\xymatrix@C=1.2em{\bar{A}\ar[r]^{-\bar{f}}&\bar{B}\ar[r]^{\bar{g}}&\bar{C}\ar[r]^{\bar{\theta}}&\Sigma
\bar{A}}$ is the right distinguished triangle given by $-f$,
obviously isomorphic to the induced triangle
$\xymatrix@C=1.5em{\bar{A}\ar[r]^{\bar{f}}&\bar{B}\ar[r]^{\bar{g}}&\bar{C}\ar[r]^{-\bar{\theta}}&\Sigma
\bar{A}}$.  \qed

In order to prove our main theorem, we also need the following
lemma.

\begin{lemma}\label{l2}
Any commutative diagram
$$\xymatrix{A\ar[r]^{f}\ar[d]^a&B\ar[r]^g\ar[d]^b&C\ar[r]^e\ar[d]^c&A[1]\ar[d]^{a[1]}\\
A'\ar[r]^{f'}&B'\ar[r]^{g'}&C'\ar[r]^{e'}&A'[1]}$$ of triangles in
$\cal T$ with $\cal X$-monomorphisms $f,f'$ of $\cal A$ induces a
commutative diagram
$$\xymatrix{\bar{A}\ar[r]^{\bar{f}}\ar[d]^{\bar{a}}&\bar{B}\ar[r]^{\bar{g}}\ar[d]^{\bar{b}}&\bar{C}\ar[r]^{\bar{h}}\ar[d]^{\bar{c}}&\Sigma \bar{A}\ar[d]^{\Sigma \bar{a}}\\
\bar{A'}\ar[r]^{\bar{f'}}&\bar{B'}\ar[r]^{\bar{g'}}&\bar{C'}\ar[r]^{\bar{h'}}&\Sigma
\bar{A'}}$$ of right induced triangles in $\cal A/[\cal X]$.
\end{lemma}

\proof We have the following three commutative diagrams
$$\xymatrix{A\ar[r]^f\ar@{=}[d]&B\ar[r]^g\ar[d]&C\ar[r]^{e}\ar[d]^{-h}&A[1]\ar@{=}[d]\\
A\ar[r]^{\alpha_A}&X_A\ar[r]^{\beta_A}&\Sigma
A\ar[r]^{\gamma_A}&A[1]}$$
$$\xymatrix{A'\ar[r]^{f'}\ar@{=}[d]&B'\ar[r]^{g'}\ar[d]&C'\ar[r]^{e'}\ar[d]^{-h'}&A'[1]\ar@{=}[d]\\
A'\ar[r]^{\alpha_A'}&X_A'\ar[r]^{\beta_A'}&\Sigma
A'\ar[r]^{\gamma_A'}&A[1]}$$
$$\xymatrix{A\ar[r]^{\alpha_A}\ar[d]^a&X_A\ar[r]^{\beta_A}\ar[d]&\Sigma A\ar[r]^{\gamma_A}\ar[d]^{a'}&A[1]\ar[d]^{a[1]}\\
A'\ar[r]^{\alpha_A'}&X_A'\ar[r]^{\beta_A'}&\Sigma
A'\ar[r]^{\gamma_A'}&A[1]'}$$ By the definition of suspension
functor $\Sigma$, we have $\Sigma \bar{a}=\bar{a'}$ and
$\gamma_A\cdot(-h)=e$, $\gamma_A\cdot(-h')=e'$. Then
$$\gamma_A'(a'h-h'c)=\gamma_A'a'h-\gamma_A'h'c=a[1]\gamma_Ah-e'c=e'c-a[1]e=0.$$
Hence $a'h-h'c$ factor through $\cal X$, i.e. $\Sigma
\bar{a}\bar{h}=\bar{h'}\bar{c}$.  \qed

\begin{theorem}
Let $\cal A$ be a full subcategory of a triangulated category $\cal
T$, $\cal X$ a covariantly finite subcategory and $\cal Y$ a
contravariantly finite subcategory of $\cal A$
\begin{itemize}
\item[\emph{(}1\emph{)}.] If the subcategory $\cal A$ has the property that any $\cal X$-monomorphism of $\cal A$
has a cone, then the subfactor category $\cal A/[\cal X]$ forms a right
triangulated category.
\item[\emph{(}2\emph{)}.] If the subcategory $\cal A$ has the property that any $\cal Y$-epimorphism of $\cal A$
has a cocone, then the subfactor category $\cal A/[\cal Y]$ forms a left
triangulated category.
\end{itemize}

\end{theorem}

\proof (\emph{1}) Let  $\nabla$ be the set of the right distinguished
triangles of $\cal A/[\cal X]$. We only need to prove the set $\nabla$ is a right
triangulation of $\cal A/[\cal X]$.
It is sufficient to check the axioms (rTR0) to (rTR4) in Definition 2.1.

(rTR0) The sextuple $\tri{\bar{A}}{\bar{A}}{0}{\Sigma
\bar{A}}{\bar{id}}{0}{0}$ is obviously a right distinguished
triangle of $\cal A/[\cal X]$ given by the identity $id:A\to A$,
since we have the following diagram of triangles.
$$\xymatrix{A\ar[d]^{id}\ar[r]^{\alpha_A}&X_A\ar[d]^{id}\ar[r]&\Sigma A\ar@{=}[d]\\
A\ar[r]^{\alpha_A}&X_A\ar[r]&\Sigma A}$$

(rTR1) Any morphism $\bar{f}:A\to B$ in $\cal A/[\cal X]$ can be
embedded into a right distinguished triangle
$\tri{\bar{A}}{\bar{B}}{\bar{C_f}}{\Sigma
\bar{A}}{\bar{f}}{\bar{g}}{\bar{h}}$ of $\cal A/[\cal X]$ since we
have the following commutative diagram of triangles in $\cal T$
$$\xymatrix{A\ar[r]^{\alpha_A}\ar[d]^f&X_A\ar[r]\ar[d]^{k}&\Sigma A\ar@{=}[d]\\
B\ar[r]^{g}&C_f\ar[r]^h&\Sigma A}$$ with $A,B,C_f,\Sigma A\in \cal
A$.

(rTR2) Let $\tri{\bar{A}}{\bar{B}}{\bar{C}_f}{\Sigma
\bar{A}}{\bar{f}}{\bar{g}}{\bar{h}}$ be a right distinguished
triangle given by $f\in \cal A$, i.e. there is the following
commutative diagram of triangles of $\cal T$.
\begin{gather}
\xymatrix{A\ar[r]^{\alpha_A}\ar[d]^f&X_A\ar[r]\ar[d]^{k}&\Sigma A\ar@{=}[d]\ar[r]^{\gamma_A}&A[1]\ar[d]^{f[1]}\\
B\ar[r]^{g}&C_f\ar[r]^h&\Sigma A\ar[r]^{f[1]\gamma_A}&B[1]}
\end{gather} Since $\alpha_A$ is $\cal X$-monic, the morphism $g\in
\cal A$ is also $\cal X$-monic by Lemma \ref{l1}. Consider the
following commutative diagram
\begin{gather}
\xymatrix{B\ar@{=}[d]\ar[r]^g&C_f\ar[d]\ar[r]^h&\Sigma A\ar[d]^u\ar[r]^{f[1]\gamma_A}&B[1]\ar@{=}[d]\\
B\ar[r]^{\alpha_B}&X_B\ar[r]^{\beta_B}&\Sigma
B\ar[r]^{\gamma_B}&B[1]}
\end{gather} of triangles in $\cal T$, i.e.
the triangle $\xymatrix{B\ar[r]^g&C_f\ar[r]^h&\Sigma
A\ar[r]^{f[1]\gamma_A}&B[1]}$ of $\cal T$ induces the right induced
triangle $\tri{\bar{B}}{\bar{C}_f}{\Sigma \bar{A}}{\Sigma
\bar{B}}{\bar{g}}{\bar{h}}{-\bar{u}}$ of $\cal A/[\cal X]$. By
composing the commutative diagrams (1) and (2), we have the
following commutative diagram
$$\xymatrix{A\ar[d]^f\ar[r]^{\alpha_A}&X_A\ar[d]\ar[r]^{\beta_A}&\Sigma A\ar[d]^u\ar[r]^{\gamma_A}&A[1]\ar[d]^{f[1]}\\
B\ar[r]^{\alpha_B}&X_B\ar[r]^{\beta_B}&\Sigma
B\ar[r]^{\gamma_B}&B[1]}$$ of triangles, then $\Sigma
\bar{f}=\bar{u}$ by the construction of suspension functor $\Sigma$
in $\cal A/[\cal X]$. Thus the right induced triangle
$\tri{\bar{B}}{\bar{C}_f}{\Sigma \bar{A}}{\Sigma
\bar{B}}{\bar{g}}{\bar{h}}{-\Sigma\bar{f}}$ belongs to $\nabla$ by
Proposition \ref{p1}.

(rTR3) Take a commutative diagram
$$\xymatrix{\bar{A}\ar[d]^{\bar{a}}\ar[r]^{\bar{f}}&\bar{B}\ar[d]^{\bar{b}}\ar[r]^{\bar{g}}&\bar{C}\ar[r]^{\bar{h}}&\Sigma \bar{A}\ar[d]^{\Sigma \bar{a}}\\
\bar{A'}\ar[r]^{\bar{f'}}&\bar{B'}\ar[r]^{\bar{g'}}&\bar{C'}\ar[r]^{\bar{h'}}&\Sigma
\bar{A'}}$$ of right distinguished triangles in $\cal A/[\cal X]$.
Owing to Proposition \ref{p1}, we can assume that these triangles
are induced, i.e. $\tri{\bar{A}}{\bar{B}}{\bar{C}}{\Sigma
\bar{A}}{\bar{f}}{\bar{g}}{\bar{h}}$ and
$\tri{\bar{A'}}{\bar{B'}}{\bar{C'}}{\Sigma
\bar{A'}}{\bar{f'}}{\bar{g'}}{\bar{h'}}$ are induced by the
triangles $\tri{A}{B}{C}{A[1]}{f}{g}{e}$ and
$\tri{A'}{B'}{C'}{A'[1]}{f'}{g'}{e'}$ with $f,f'\in \cal A$ being
$\cal X$-monic, respectively. By the construction of right induced
triangles, we have the following (not necessarily commutative)
diagram
$$\xymatrix{A\ar[d]^{a}\ar[r]^f&B\ar[d]^{b}\ar[r]^g&C\ar[r]^e&A[1]\ar[d]^{a[1]}\\
A'\ar[r]^{f'}&B'\ar[r]^{g'}&C'\ar[r]^{e'}&A'[1]}$$ of triangles in
$\cal T$. Since $\bar{b}\bar{f}=\bar{f}'\bar{a}$ holds, $bf-f'a$
factors through $\cal X$, i.e. there exist morphisms $t:A\to X$ and
$t':X\to Y$ for some $X\in \cal X$ such that $bf-f'a=t't$. Because
$f$ is $\cal X$-monic, there exists a morphism $s:B\to X$ such that
$t=sf$. Then $bf-f'a=(t's)f$ with $t's\in [\cal X](B,B')$. Set
$b'=b-t's$, then we have $b'f=bf-t'sf=f'a$. Then we get the
following commutative diagram
$$\xymatrix{A\ar[d]^{a}\ar[r]^f&B\ar[d]^{b'}\ar[r]^g&C\ar@{-->}[d]^{\exists c}\ar[r]^e&A[1]\ar[d]^{a[1]}\\
A'\ar[r]^{f'}&B'\ar[r]^{g'}&C'\ar[r]^{e'}&A'[1]}$$ of triangles in
$\cal T$ with $f,f'\in \cal A$ being $\cal X$-monic. Thus the
assertion follows from Lemma \ref{l2}, since $\bar{b}'=\bar{b}$.

(rTR4) As the argumentation in (rTR3), we can assume by Proposition
\ref{p1}, the triangles are induced, i.e.
$\rtri{\bar{A}}{\bar{B}}{\bar{C}}{\bar{f}}{\bar{g}}{\bar{h}}$,
$\rtri{\bar{B}}{\bar{X}}{\bar{Y}}{\bar{a}}{\bar{b}}{\bar{c}}$ and
$\rtri{\bar{A}}{\bar{X}}{\bar{Z}}{\bar{a}\bar{f}}{\bar{d}}{\bar{e}}$
are induced by the triangles $\tri{A}{B}{C}{A[1]}{f}{g}{h'}$,\\
$\tri{B}{X}{Y}{B[1]}{a}{b}{c'}$ and $\tri{A}{X}{Z}{A[1]}{af}{d}{e'}$
with $f,a,af\in \cal A$ being $\cal X$-monic. Then by the octahedral
axioms in $\cal T$, we have the following commutative diagrams
$$\xymatrix{A\ar[r]^f\ar@{=}[d]&B\ar[r]^g\ar[d]^a&C\ar[r]^{h'}\ar[d]^s&A[1]\ar@{=}[d]\\
A\ar[r]^{af}&X\ar[r]^d\ar[d]^b&Z\ar[r]^{e'}\ar[d]^{t}&A[1]\\
&Y\ar@{=}[r]\ar[d]^{c'}&Y\ar[d]\\
&B[1]\ar[r]&C[1]}$$
$$\xymatrix{A\ar[r]^{af}\ar[d]^f&X\ar[r]^d\ar@{=}[d]&Z\ar[r]^{e'}\ar[d]^t&A[1]\ar[d]^{f[1]}\\
B\ar[r]^a&X\ar[r]^b&Y\ar[r]^{c'}&B[1]}$$ of triangles in $\cal T$
with $s$ being $\cal X$-monic by Lemma \ref{l1}. Thus by Lemma
\ref{l2}, we have the following commutative diagrams
$$\xymatrix{\bar{A}\ar[r]^{\bar{f}}\ar@{=}[d]&\bar{B}\ar[r]^{\bar{g}}\ar[d]^{\bar{a}}&\bar{C}\ar[r]^{\bar{h}}\ar[d]^{\bar{s}}&\Sigma \bar{A}\ar@{=}[d]\\
\bar{A}\ar[r]^{\bar{a}\bar{f}}&\bar{X}\ar[r]^{\bar{d}}\ar[d]^{\bar{b}}&\bar{Z}\ar[r]^{\bar{e}}\ar[d]^{\bar{t}}&\Sigma \bar{A}\\
&\bar{Y}\ar@{=}[r]\ar[d]^{\bar{c}}&\bar{Y}\ar[d]\\
&\Sigma \bar{B}\ar[r]&\Sigma \bar{C}}$$
$$\xymatrix{\bar{A}\ar[r]^{\bar{a}\bar{f}}\ar[d]^{\bar{f}}&\bar{X}\ar[r]^{\bar{d}}\ar@{=}[d]&\bar{Z}\ar[r]^{\bar{e}}\ar[d]^{\bar{t}}&\Sigma \bar{A}\ar[d]^{\Sigma \bar{f}}\\
\bar{B}\ar[r]^{\bar{a}}&\bar{X}\ar[r]^{\bar{b}}&\bar{Y}\ar[r]^{\bar{c}}&\Sigma
\bar{B}}$$ of right triangles in $\cal A/[\cal X]$.

(\emph{2}) Dual to (1). Let $\triangle$ be the set of the left distinguished
triangles of $\cal A/[\cal Y]$. We only need to prove the set $\triangle$ is a left
triangulation of $\cal A/[\cal Y]$.
 \qed

\begin{corollary}Let $\cal A$ be a full subcategory of a triangulated category $\cal
T$, $\cal X$ a functorially finite subcategory of $\cal A$. If $\cal
A$ has the properties that any $\cal X$-monomorphism of $\cal A$ has
a cone and any $\cal X$-epimorphism has a cocone. Then the subfactor
category $\cal A/[\cal X]$ admits a pretriangulated structure in the
sense of \cite{br}.
\end{corollary}

\proof An easy modification of results of \cite{br,b} in our
setting. \qed

\begin{remark}
A pretriangulation $(\Omega,\Sigma,\triangle,\nabla)$ of $\cal C$
becomes a triangulation if and only if $\triangle=\nabla$ and
$\Omega=\Sigma^{-1}$.
\end{remark}

\begin{example} \label{e1}Let $\cal T$ be a triangulated category and $\cal X\subset\cal
A$ be subcategories of $\cal T$. If $\cal A$ and $\cal X$ satisfy
the conditions of \cite{iy}:

(Z1) $\cal A$ is extension closed,

(Z2) $ (\cal A,\cal A)$ forms an $\cal X$-mutation pair. \\
Then by the definition of $\cal X$-mutation in \cite{iy}, $\cal X$
is a functorially finite subcategory of $\cal A$ and $\cal A$ has
the properties that any $\cal X$-monomorphism has a cone and any
$\cal X$-epimorphism has a cocone by Lemma 4.3.(2) in \cite{iy} and
its dual. Moreover, the suspension functor coincides with the loop
functor and is an auto-equivalence of $\cal A/[\cal X]$ [see
Proposition 2.6 \cite{iy}. Hence the subfactor category $\cal
A/[\cal X]$ forms a triangulated category, that is the Theorem 4.2
of \cite{iy}.
\end{example}
In addition, we can prove the converse is also true.
\begin{theorem}Let $\cal T$ be a triangulated category and $\cal X\subset\cal
A$ be subcategories of $\cal T$. $\cal A$ and $\cal X$ satisfy the
conditions (Z1) and (Z2) if and only if the subfactor category $\cal
A/[\cal X]$ forms a triangulated category where $\cal X$ is a
functorially finite subcategory of $\cal A$ satisfying $\cal
{(X,X}[1])=0$ and $\cal A$ has the properties that any $\cal
X$-monomorphism has a cone and any $\cal X$-epimorphism has a
cocone.
\end{theorem}

\proof $``\Rightarrow"$ Example \ref{e1}.

$``\Leftarrow"$ First of all, we claim $(\cal X, \cal A[1])=0$.
Indeed, since the pretriangulated category $\cal A/[\cal X]$ forms a
triangulated category, then $\triangle=\nabla$ and
$\Omega=\Sigma^{-1}$ is an auto-equivalence of $\cal A/[\cal X]$.
For any object $A$ in $\cal A$, we have the following commutative
diagram
$$\xymatrix{A\ar[d]\ar[r]^{\alpha}&X_A\ar[d]\ar[r]&\Sigma A\ar@{=}[d]\ar[r]&A[1]\ar[d]\\
\Omega\Sigma A\ar[r]&X_{\Sigma A}\ar[r]^{\beta}&\Sigma
A\ar[r]&\Omega\Sigma A[1]}$$ of triangles with $\alpha$ a left $\cal
X$-approximation of $A$ and $\beta$ a right $\cal X$-approximation
of $\Sigma A$. Since $\beta$ is a right $\cal X$-approximation and
\mbox{$(\cal X,\cal X[1])=0$, we have $(\cal X,\Omega\Sigma A[1])=0$}.
Then by $\Omega\Sigma \bar{A}\cong \bar{A}$, i.e. $\Omega\Sigma
A\oplus X'=A\oplus X''$, we get $(\cal X,A[1])=0$. Hence $(\cal X,
\cal A[1])=0$ holds. Dually we can get $(\cal A, \cal X[1])=0$.

Let $\tri{A_2}{B}{A_1}{A_2[1]}{}{}{}$ be a triangle of $\cal T$ with
$A_1,A_2 \in \cal A$. Noting that $(\cal X, \cal A[1])=0$, we have
the following diagram
$$\xymatrix{&A_3\ar@{-->}[d]^s\ar@{=}[r]&A_3\ar[d]^{\beta}\\
A_2\ar@{=}[d]\ar[r]^{\binom{1}{0}\quad}&A_2\oplus X\ar@{-->}[d]^t\ar[r]^{\quad(0,1)}&X\ar[d]^{\alpha}\ar[r]^0&A_2[1]\ar@{=}[d]\\
A_2\ar[r]&B\ar[r]&A_1\ar[r]&A_2[1]}$$ of triangles in $\cal T$ where
$\alpha$ is a right $\cal X$-approximation of $A_1$. Easy
computation allow us to have $s=\binom{l}{\beta}$ for some $l$.
Since $\alpha$ is a right $\cal X$-approximation and $(\cal A, \cal
X[1])=0$, we have $A_3\in \cal A$ and $\beta$ is a left $\cal
X$-approximation of $A_3$. By Lemma \ref{l1}, $s=\binom{l}{\beta}$
is $\cal X$-monic with $A_3,A_2,X\in \cal A$, hence $B\in \cal A$,
i.e. $\cal A$ is closed under extensions.

Since $\Sigma$ is an auto-equivalence of $\cal A/[\cal X]$, for
$\forall C\in \cal A$, there exists an object $A$ such that $\Sigma
\bar{A}=\bar{C}$, i.e. satisfying the following triangle
$$\tri{A}{X_A}{\Sigma A}{A[1]}{\alpha_A}{}{}$$ with $\alpha_A$ a $\cal
X$-approximation of $A$ and $X_A\in \cal X$. Then $C\oplus
X_1=\Sigma A\oplus X_2 $ holds for some $X_1,X_2\in \cal X$. Since
$\Sigma A\oplus X_2\in \mu^{-1}(\cal A;\cal X)$, $(\cal X, \cal
A[1])=0$, then $\mu^{-1}(\cal A;\cal X)=(\cal X*\cal A[1])\cap
^\bot\cal X[1]$ is closed under direct summands by Proposition 2.1
of \cite{iy}. Thus $C\in \mu^{-1}(\cal A;\cal X)$, i.e. $\cal
X\subset\cal A\subset \mu^{-1}(\cal A;\cal X)$.

Dually we can also get $\cal X\subset \cal A \subset \mu(\cal A;\cal
X)$. Hence $(\cal A,\cal A)$ forms a $\cal X$-mutation pair.  \qed

\begin{comment}
\section*{Acknowledgement}
\end{comment}

\end{document}